\documentclass{amsart}
\usepackage{amsmath,amssymb}
\usepackage{yhmath}

\newcommand{\bdis}{\begin{displaymath}}
\newcommand{\edis}{\end{displaymath}}
\newcommand{\be}{\begin{equation}}
\newcommand{\ee}{\end{equation}}
\newcommand{\mbb}{\mathbb}
\newcommand{\mcal}{\mathcal}

\newcommand{\vp}{\varphi}

\newcommand{\zf}{\zeta\left(\frac{1}{2}+it\right)}

\newcommand{\zfvpr}{\zeta\left(\frac{1}{2}+i\vp_1^r(t)\right)}

\newcommand{\FR}{\frac{x^n+y^n}{z^n}}

\DeclareMathOperator{\li}{li}

\theoremstyle{definition}

\theoremstyle{remark}
\newtheorem{remark}[]{Remark}

\newtheorem*{mydef11}{{\bf Theorem 1}}

\newtheorem*{mydef12}{{\bf Theorem 2}}

\newtheorem*{mydef13}{{\bf Theorem 3}}

\newtheorem*{mydef14}{{\bf Theorem 4}}

\newtheorem*{mydef15}{{\bf Theorem 5}} 

\newtheorem*{mydef16}{{\bf Theorem 6}} 

\newtheorem*{mydef17}{{\bf Theorem 7}} 

\newtheorem*{mydef18}{{\bf Theorem 8}} 

\newtheorem*{mydef19}{{\bf Theorem 9}}

\newtheorem*{mydef41}{{\bf Corollary 1}}

\newtheorem*{mydef42}{{\bf Corollary 2}}

\newtheorem*{mydef51}{{\bf Lemma 1}}

\newtheorem*{mydef81}{{\bf Property 1}}

\newtheorem*{mydef82}{{\bf Property 2}}

\numberwithin{equation}{section}

\begin{document}

\title[Jacob's ladders, point of contact of the remainder in the prime-number law \dots]{Jacob's ladders, point of contact of the remainder in the prime-number law with the Fermat-Wiles theorem and  multiplicative puzzles on some sets of integrals}

\author{Jan Moser}

\address{Department of Mathematical Analysis and Numerical Mathematics, Comenius University, Mlynska Dolina M105, 842 48 Bratislava, SLOVAKIA}

\email{jan.mozer@fmph.uniba.sk}


\begin{abstract}
In this paper we prove, on the Riemann hypothesis, the existence of such increments of the Ingham integral (1932) that generate new functionals together with corresponding new $P\zeta$-equivalents of the Fermat-Wiles theorem. We obtain also new results in this direction. 
\end{abstract}
\maketitle

\section{Introduction} 

\subsection{} 

Let us remind that 
\be \label{1.1} 
\pi(t)=P(t)+\int_0^t\frac{{\rm d}t}{\ln t}=\li(t)+P(t), 
\ee  
where $\pi(t)$ is the prime-counting function. Littlewood proved the following fundamental result in 1914: the remainder $P(t)$,\ $t\geq 2$ changes the sign infinitely many times, see \cite{4}. 

However, Ingham proved in 1932 that, on the Riemann hypothesis, the inequality 
\be \label{1.2} 
\int_2^YP(t){\rm d}t<0 
\ee  
holds true for every sufficiently big $Y$, see \cite{3}, p. 106, i.e. the remainder $P(t)$ is negative in mean. 

\subsection{} 

In the direction of the Ingham's integral in (\ref{1.2}) we have proved, see \cite{7}, the following theorem. On the Riemann hypothesis there is such a positive function $N(\delta)$, $\delta\in (0,\Delta)$ that 
\be \label{1.3} 
\int_2^{N(\delta)}P(t){\rm d}t=-\frac{a}{\delta}+\mcal{O}(1),\ \delta\in(0,\Delta),\ a=e^{4.5}\approx 90.017, 
\ee 
where $\mcal{O}(1)$ stands for a continuous and bounded function. 

\begin{remark}
Since, see (\ref{1.3}), 
\bdis 
\lim_{\delta\to 0^+}\int_2^{N(\delta)}\{-P(t)\}{\rm d}t=+\infty, 
\edis  
then 
\be \label{1.4} 
\lim_{\delta\to 0^+}N(\delta)=+\infty. 
\ee 
\end{remark} 

\subsection{} 

In this paper we prove, on the Riemann hypothesis, the existence of such increments of the Ingham's integral in (\ref{1.2}) which generate new functional together with corresponding $P\zeta$-equivalents of Fermat-Wiles theorem and also some other results. 

For example, we obtain the following functional 
\be \label{1.5} 
\begin{split}
& \lim_{\tau\to\infty}\frac{1}{\tau}\left\{ 
\int_{N_1(x\tau)}^{N_1([x\tau]^1)}\{-P(t)\}{\rm d}t\times \int_{[x\tau]^1}^{[x\tau+\frac{1}{a(1-c)}]^1}\left|\zf\right|^2{\rm d}t
\right\}=x
\end{split}
\ee 
for every fixed $x>0$, where 
\be \label{1.6} 
N_1(G)=N\left(\frac{1}{G}\right),\ [G]^1=\vp_1^{-1}(G). 
\ee 
In the special case 
\be \label{1.7} 
x\to \FR,\ x,y,z,n\in\mbb{N},\ n\geq 3
\ee 
of the Fermat's rationals we have the following result: on the Riemann hypothesis the following $P\zeta$-condition 
\be \label{1.8} 
\begin{split}
& \lim_{\tau\to\infty}\frac{1}{\tau}\left\{ 
\int_{N_1(\FR\tau)}^{N_1([\FR\tau]^1)}\{-P(t)\}{\rm d}t\times \right. \\ 
& \left. \int_{[\FR\tau]^1}^{[\FR\tau+\frac{1}{a(1-c)}]^1}\left|\zf\right|^2{\rm d}t
\right\}\not=1 
\end{split}
\ee 
on the set of all Fermat's rationals gives new $P\zeta$-equivalent of the Fermat-Wiles theorem. 

Next, we have obtained the asymptotic factorization formula of this sort: 
\be \label{1.9} 
\begin{split}
& \int_T^{[T]^1}\left|\zf\right|^2{\rm d}t\sim \\ 
& \int_{[T]^1}^{[T+\frac{1}{a}]^1}\left|\zf\right|^2{\rm d}t\times\int_{N_1(T)}^{N_1([T]^1)}\{-P(t)\}{\rm d}t,\ T\to\infty. 
\end{split}
\ee 

\begin{remark}
Let us remind explicitly that the $P\zeta$-condition represents the 5-fold point of contact namely between the Riemann zeta-function, the Riemann hypothesis, the prime-number law, the Jacob's ladders and the Fermat-Wiles theorem. 
\end{remark} 

\subsection{} 

On the Riemann hypothesis we have obtained also the following results: 
\begin{itemize}
	\item[(A)] the asymptotic mapping 
	\bdis 
	-P(t)\xrightarrow{\mcal{T}}\pi(t) 
	\edis 
	defined by the formula 
	\be \label{1.10} 
	\int_{N_1(\frac{\tau}{a(1-c)})}^{N_1([\frac{\tau}{a(1-c)})]^1)}\{-P(t)\}{\rm d}t\sim\pi(\tau),\ \tau\to\infty, 
	\ee 
	\item[(B)] the new double integral\footnote{Comp. (\ref{1.9}).} 
	\be \label{1.11} 
	\begin{split}
	& \iint_{Q_2(T)}\{-P(u)\}\left|\zeta\left(\frac{1}{2}+iv\right)\right|^2{\rm d}u{\rm d}v\sim(1-c)T,\ T\to\infty, \\ 
	& Q_2(T)=\{(u,v):\ u\in (N_1(T),N_1([T]^1)),\ v\in ([T]^1,[T+\frac{1}{a}]^1)\}, 
	\end{split} 
	\ee 
	that expresses the result of new type of interaction 
	\bdis 
	-P(u)\leftrightarrow \left|\zeta\left(\frac{1}{2}+iv\right)\right|^2, 
	\edis  
	\item[(C)] and, for example, the following multiplicative puzzle 
	\be \label{1.12} 
	\begin{split}
	& \int_2^{N_2^4(\ln^4T)}\{-P(t)\}{\rm d}t\sim \\ 
	& \int_{[T]^1}^{[T+1]^1}\left|\zf\right|^2{\rm d}t\times \int_{[T]^1}^{[T+e]^1}\left|\zf\right|^2{\rm d}t\times \\ 
	& \int_{[T]^1}^{[T+e^{1.5}]^1}\left|\zf\right|^2{\rm d}t\times \int_{[T]^1}^{[T+e^2]^1}\left|\zf\right|^2{\rm d}t
	\end{split}
	\ee 
	on the set of corresponding integrals. 
\end{itemize}

\section{Jacob's ladders: notions and basic geometrical properties}  

\subsection{}

In this paper we use the following notions of our works \cite{5}, \cite{6},  \cite{8} -- \cite{10}: 
\begin{itemize}
\item[{\tt (a)}] Jacob's ladder $\vp_1(T)$, 
\item[{\tt (b)}] direct iterations of Jacob's ladders 
\bdis 
\begin{split}
	& \vp_1^0(t)=t,\ \vp_1^1(t)=\vp_1(t),\ \vp_1^2(t)=\vp_1(\vp_1(t)),\dots , \\ 
	& \vp_1^k(t)=\vp_1(\vp_1^{k-1}(t))
\end{split}
\edis 
for every fixed natural number $k$, 
\item[{\tt (c)}] reverse iterations of Jacob's ladders 
\be \label{2.1}  
\begin{split}
	& \vp_1^{-1}(T)=\overset{1}{T},\ \vp_1^{-2}(T)=\vp_1^{-1}(\overset{1}{T})=\overset{2}{T},\dots, \\ 
	& \vp_1^{-r}(T)=\vp_1^{-1}(\overset{r-1}{T})=\overset{r}{T},\ r=1,\dots,k, 
\end{split} 
\ee   
where, for example, 
\be \label{2.2} 
\vp_1(\overset{r}{T})=\overset{r-1}{T}
\ee  
for every fixed $k\in\mbb{N}$ and every sufficiently big $T>0$. We also use the properties of the reverse iterations listed below.  
\be \label{2.3}
\overset{r}{T}-\overset{r-1}{T}\sim(1-c)\pi(\overset{r}{T});\ \pi(\overset{r}{T})\sim\frac{\overset{r}{T}}{\ln \overset{r}{T}},\ r=1,\dots,k,\ T\to\infty,  
\ee 
\be \label{2.4} 
\overset{0}{T}=T<\overset{1}{T}(T)<\overset{2}{T}(T)<\dots<\overset{k}{T}(T), 
\ee 
and 
\be \label{2.5} 
T\sim \overset{1}{T}\sim \overset{2}{T}\sim \dots\sim \overset{k}{T},\ T\to\infty.   
\ee  
\end{itemize} 

\begin{remark}
	The asymptotic behaviour of the points 
	\bdis 
	\{T,\overset{1}{T},\dots,\overset{k}{T}\}
	\edis  
	is as follows: at $T\to\infty$ these points recede unboundedly each from other and all together are receding to infinity. Hence, the set of these points behaves at $T\to\infty$ as one-dimensional Friedmann-Hubble expanding Universe. 
\end{remark}  

\subsection{} 

Let us remind that we have proved\footnote{See \cite{10}, (3.4).} the existence of almost linear increments 
\be \label{2.6} 
\begin{split}
& \int_{\overset{r-1}{T}}^{\overset{r}{T}}\left|\zf\right|^2{\rm d}t\sim (1-c)\overset{r-1}{T}, \\ 
& r=1,\dots,k,\ T\to\infty,\ \overset{r}{T}=\overset{r}{T}(T)=\vp_1^{-r}(T)
\end{split} 
\ee 
for the Hardy-Littlewood integral (1918) 
\be \label{2.7} 
J(T)=\int_0^T\left|\zf\right|^2{\rm d}t. 
\ee  

For completeness, we give here some basic geometrical properties related to Jacob's ladders. These are generated by the sequence 
\be \label{2.8} 
T\to \left\{\overset{r}{T}(T)\right\}_{r=1}^k
\ee 
of reverse iterations of the Jacob's ladders for every sufficiently big $T>0$ and every fixed $k\in\mbb{N}$. 

\begin{mydef81}
The sequence (\ref{2.8}) defines a partition of the segment $[T,\overset{k}{T}]$ as follows 
\be \label{2.9} 
|[T,\overset{k}{T}]|=\sum_{r=1}^k|[\overset{r-1}{T},\overset{r}{T}]|
\ee 
on the asymptotically equidistant parts 
\be \label{2.10} 
\begin{split}
& \overset{r}{T}-\overset{r-1}{T}\sim \overset{r+1}{T}-\overset{r}{T}, \\ 
& r=1,\dots,k-1,\ T\to\infty. 
\end{split}
\ee 
\end{mydef81} 

\begin{mydef82}
Simultaneously with the Property 1, the sequence (\ref{2.8}) defines the partition of the integral 
\be \label{2.11} 
\int_T^{\overset{k}{T}}\left|\zf\right|^2{\rm d}t
\ee 
into the parts 
\be \label{2.12} 
\int_T^{\overset{k}{T}}\left|\zf\right|^2{\rm d}t=\sum_{r=1}^k\int_{\overset{r-1}{T}}^{\overset{r}{T}}\left|\zf\right|^2{\rm d}t, 
\ee 
that are asymptotically equal 
\be \label{2.13} 
\int_{\overset{r-1}{T}}^{\overset{r}{T}}\left|\zf\right|^2{\rm d}t\sim \int_{\overset{r}{T}}^{\overset{r+1}{T}}\left|\zf\right|^2{\rm d}t,\ T\to\infty. 
\ee 
\end{mydef82} 

It is clear, that (\ref{2.10}) follows from (\ref{2.3}) and (\ref{2.5}) since 
\be \label{2.14} 
\overset{r}{T}-\overset{r-1}{T}\sim (1-c)\frac{\overset{r}{T}}{\ln \overset{r}{T}}\sim (1-c)\frac{T}{\ln T},\ r=1,\dots,k, 
\ee  
while our eq. (\ref{2.13}) follows from (\ref{2.6}) and (\ref{2.5}).  

\section{The first $P\zeta$-functional and corresponding $P\zeta$-equivalent generated by the formula (\ref{1.3})} 

\subsection{} 

If we put in (\ref{1.3})\footnote{See (\ref{2.4}), (\ref{2.5}).} consecutively, 
\be \label{3.1} 
\delta=\frac{1}{\overset{r-1}{T}(T)},\ N\left(\frac{1}{\overset{r-1}{T}}\right)=N_1(\overset{r-1}{T}),\ r=1,2,\dots,k+1,\ \overset{0}{T}=T, 
\ee  
where, of course, 
\be \label{3.2} 
\delta\to 0^+ \ \Leftrightarrow \ \overset{r-1}{T}(T)\to +\infty, 
\ee  
then we obtain, for example, 
\be \label{3.3} 
\begin{split}
& \int_2^{N_1(\overset{r-1}{T})}\{-P(t)\}{\rm d}t=a\overset{r-1}{T}+\mcal{O}(1), \\ 
& \int_2^{N_1(\overset{r}{T})}\{-P(t)\}{\rm d}t=a\overset{r}{T}+\mcal{O}(1), 
\end{split} 
\ee 
i.e. we have obtained the following increments of the integral (\ref{1.3}). 
\be \label{3.4} 
\int_{N_1(\overset{r-1}{T})}^{N_1(\overset{r}{T})}\{-P(t)\}{\rm d}t=a(\overset{r}{T}-\overset{r-1}{T})+\mcal{O}(1), \ r=1,\dots,k,\ T\to\infty. 
\ee 
Now, we have by (\ref{2.14}) the following result. 

\begin{mydef51} 
	On the Riemann hypothesis 
\be \label{3.5} 
\int_{{N_1(\overset{r-1}{T})}}^{N_1(\overset{r}{T})}\{-P(t)\}{\rm d}t=\{1+o(1)\}a(1-c)\frac{T}{\ln T},\ r=1,\dots,k,\ T\to\infty. 
\ee 
\end{mydef51} 
Consequently, we have the following geometrical property from (\ref{3.5}). 

\begin{mydef11}
On the Riemann hypothesis and for every sufficiently big $T>0$ and for every fixed $k\in\mbb{N}$ the sequence 
\be \label{3.6} 
\{N_1(\overset{r}{T})\}_{r=1}^k:\ T\to \{\overset{r}{T}(T)\}_{r=1}^k\to \{N_1(\overset{r}{T}(T))\}_{r=1}^k
\ee 
defines the partition of the integral 
\be \label{3.7} 
\int_{{N_1(T)}}^{N_1(\overset{k}{T})}\{-P(t)\}{\rm d}t
\ee 
into the parts 
\be \label{3.8} 
\int_{{N_1(T)}}^{N_1(\overset{k}{T})}\{-P(t)\}{\rm d}=\sum_{r=1}^k\int_{{N_1(\overset{r-1}{T})}}^{N_1(\overset{r}{T})}\{-P(t)\}{\rm d}t
\ee 
that are asymptotically equal 
\be \label{3.9} 
\begin{split} 
& \int_{{N_1(\overset{r-1}{T})}}^{N_1(\overset{r}{T})}\{-P(t)\}{\rm d}t\sim \int_{{N_1(\overset{r}{T})}}^{N_1(\overset{r+1}{T})}\{-P(t)\}{\rm d}t, \\ 
& r=1,\dots,k-1,\ T\to\infty. 
\end{split} 
\ee 
\end{mydef11} 

\begin{remark}
It is clear that this Theorem is an analogue of the Property 2, see (\ref{2.11}) -- (\ref{2.13}). 
\end{remark} 

\subsection{} 

Next, we use, as the basis, the formula (\ref{3.5}) with $r=1$ in the form 
\be \label{3.10} 
\int_{{N_1(T)}}^{N_1(\overset{1}{T})}\{-P(t)\}{\rm d}t=\{1+o(1)\}a(1-c)\frac{T}{\ln T},\ T\to\infty . 
\ee  
Further, let us remind our formula\footnote{See \cite{11}, (3.18), $f_m\equiv 1$.} 
\be \label{3.11} 
\begin{split}
& \int_{[T]^k}^{[T+2l]^k}\prod_{r=0}^{k-1}\left|\zfvpr\right|^2{\rm d}t= 
 \left\{1+\mcal{O}\left(\frac{\ln \ln T}{\ln T}\right)\right\}2l\ln^kT
\end{split}
\ee 
for every fixed $l>0$ and $k\in\mbb{N}$. Now, we use this formula for $k=1$ in the form 
\be \label{3.12} 
2l\ln T=\{1+o(1)\}\int_{[T]^1}^{[T+2l]^1}\left|\zf\right|^2{\rm d}t
\ee 
in (\ref{3.10}) that gives us the following result 
\be \label{3.13} 
\begin{split}
& \int_{{N_1(T)}}^{N_1([T]^1)}\{-P(t)\}{\rm d}t \times \int_{[T]^1}^{[T+2l]^1}\left|\zf\right|^2{\rm d}t= \\ 
& \{1+o(1)\}2la(1-c)T,\ T\to\infty. 
\end{split}
\ee 
Consequently, by making the choice for $l$: 
\be \label{3.14} 
l=\bar{l}=\frac{1}{2a(1-c)}, 
\ee  
and then, by the substitution 
\be \label{3.15} 
T=x\tau,\ x>0 
\ee 
in the formula (\ref{3.13}), we obtain the following functional. 

\begin{mydef12} 
On the Riemann hypothesis it is true that 
\be \label{3.16} 
\begin{split}
& \lim_{\tau\to\infty}\frac{1}{\tau}\left\{
\int_{{N_1(x\tau)}}^{N_1([x\tau]^1)}\{-P(t)\}{\rm d}t\times \int_{[x\tau]^1}^{[x\tau+\frac{1}{a(1-c)}]^1}\left|\zf\right|^2{\rm d}t
\right\}=x
\end{split}
\ee 
for every fixed $x>0$ . 
\end{mydef12}  

\subsection{} 

In the special case of the Fermat's rationals, see (\ref{1.7}), we obtain the following statement. 

\begin{mydef41}
\be \label{3.17} 
\begin{split}
& \lim_{\tau\to\infty}\frac{1}{\tau}\left\{
\int_{{N_1(\FR\tau)}}^{N_1([\FR\tau]^1)}\{-P(t)\}{\rm d}t\times \right. \\
& \left. \int_{[\FR\tau]^1}^{[\FR\tau+\frac{1}{a(1-c)}]^1}\left|\zf\right|^2{\rm d}t
\right\} = \FR
\end{split}
\ee 
for every fixed Fermat's rational. 
\end{mydef41} 

Next Theorem follows from (\ref{3.17}). 

\begin{mydef13}
On the Riemann hypothesis, the following $P\zeta$-condition  
\be \label{3.18} 
\begin{split}
	& \lim_{\tau\to\infty}\frac{1}{\tau}\left\{
	\int_{{N_1(\FR\tau)}}^{N_1([\FR\tau]^1)}\{-P(t)\}{\rm d}t\times \right. \\
	& \left. \int_{[\FR\tau]^1}^{[\FR\tau+\frac{1}{a(1-c)}]^1}\left|\zf\right|^2{\rm d}t
	\right\} \not=1
\end{split}
\ee 
on the set of all Fermat's rationals gives the new $P\zeta$-equivalent of the Fermat-Wiles theorem. 
\end{mydef13} 

\section{Riemann hypothesis, transformation $\{-P(t)\}\to\pi(t)$ and new $P\zeta$-equivalent of the Fermat-Wiles theorem} 

\subsection{} 

Now, we use the substitution 
\be \label{4.1} 
a(1-c)T=\tau 
\ee 
in the formula (\ref{3.10}). Since 
\be \label{4.2} 
\ln T=\left\{1+\mcal{O}\left(\frac{1}{\ln \tau}\right)\right\}\ln\tau, 
\ee 
we have that 
\be \label{4.3} 
\begin{split}
& \int_{N_1(\frac{\tau}{a(1-c)})}^{N_1([\frac{\tau}{a(1-c)}]^1)}\{-P(t)\}{\rm d}t\sim \frac{\tau}{\ln\tau},\ \tau\to\infty. 
\end{split}
\ee 
Next, the prime-number law implies that 
\be \label{4.4} 
\frac{\tau}{\ln\tau}=\{1+o(1)\}\pi(\tau),\ \tau\to\infty. 
\ee 
Consequently the next Theorem follows from (\ref{4.3}) by (\ref{4.4}).  

\begin{mydef14} 
On the Riemann hypothesis the following asymptotic formula holds true 
\be \label{4.5} 
\int_{N_1(\frac{\tau}{a(1-c)})}^{N_1([\frac{\tau}{a(1-c)}]^1)}\{-P(t)\}{\rm d}t\sim \pi(t),\ \tau\to\infty. 
\ee 
\end{mydef14} 

\begin{remark}
That is, on the Riemann hypothesis, the parametric integral (\ref{4.5}) defines the the asymptotic mapping 
\be \label{4.6} 
-P(t)\xrightarrow{\mcal{T}}\pi(t), 
\ee  
i.e. the asymptotic mapping of the minus remainder onto the prime-counting function, comp. (\ref{1.1}). 
\end{remark} 

\begin{remark}
Of course, our parametric integral in (\ref{4.5}) represents the continuum set of the corresponding increments of the Ingham's integral in (\ref{1.2}). 
\end{remark} 

\subsection{} 

Next, we use the formula (\ref{4.3}) in the form 
\be \label{4.7} 
\int_{N_1(\frac{\tau}{a(1-c)})}^{N_1([\frac{\tau}{a(1-c)}]^1)}\{-P(t)\}{\rm d}t=\{1+o(1)\}\frac{\tau}{\ln\tau},\ \tau\to\infty, 
\ee  
and the substitution 
\be \label{4.8} 
\tau=x\rho,\ x>0;\ \{\tau\to+\infty\} \ \Leftrightarrow \ \{\rho\to+\infty\}
\ee 
in (\ref{4.7}). Since 
\be \label{4.9} 
\ln\tau=\left\{1+\frac{\ln x}{\ln\rho}\right\}\ln\rho=\{1+o(1)\}\ln\rho, 
\ee  
for every fixed $x>0$, then it follows from (\ref{4.7}): 
\be \label{4.10} 
\int_{N_1(\frac{x\rho}{a(1-c)})}^{N_1([\frac{x\rho}{a(1-c)}]^1)}\{-P(t)\}{\rm d}t=\{1+o(1)\}\frac{x\rho}{\ln\rho},\ \rho\to\infty. 
\ee 
Now, we use (\ref{4.4}) with $\tau\to\rho$ and obtain the following functional.  

\begin{mydef15}
On the Riemann hypothesis it is true that 
\be \label{4.11} 
\lim_{\rho\to\infty}\frac{1}{\pi(\rho)}\int_{N_1(\frac{x\rho}{a(1-c)})}^{N_1([\frac{x\rho}{a(1-c)}]^1)}\{-P(t)\}{\rm d}t=x 
\ee  
for every fixed $x>0$. 
\end{mydef15} 

\subsection{} 

In the special case of Fermat's rationals, see (\ref{1.7}), we obtain the following result\footnote{Comp. (\ref{3.16}) -- (\ref{3.18}).}. 

\begin{mydef16}
On the Riemann hypothesis the following $P\zeta$-condition 
\be\label{4.12} 
\lim_{\rho\to\infty}\frac{1}{\pi(\rho)}\int_{N_1(\frac{\rho}{a(1-c)}\FR)}^{N_1([\frac{\rho}{a(1-c)}\FR]^1)}\{-P(t)\}{\rm d}t\not=1 
\ee  
on the set of all Fermat's rationals gives the new $P\zeta$-equivalent of the Fermat-Wiles theorem. 
\end{mydef16} 

\begin{remark}
Namely, $P\zeta$-equivalent, since 
\bdis 
[G]^1=\vp_1^{-1}(G), 
\edis  
and the Jacob's ladder was generated by the Riemann $\zf$-function, see \cite{5}, and also by the Hardy-Littlewood integral (1918), see \cite{2}. 
\end{remark} 

\section{Some factorization formulas based on the Riemann hypothesis and the Ingham's integral} 

\subsection{} 

Now we use our two formulae (\ref{3.12}) and (\ref{2.5}) with $r=1$ in\footnote{See (\ref{3.10}).} 
\bdis 
\int_{{N_1(T)}}^{N_1([T]^1)}\{-P(t)\}{\rm d}t = \{1+o(1)\}a(1-c)\frac{T}{\ln T}, 
\edis  
and this gives 
\be \label{5.1} 
\begin{split}
& \int_{{N_1(T)}}^{N_1([T]^1)}\{-P(t)\}{\rm d}t\times \int_{[T]^1}^{[T+2l]^1}\left|\zf\right|^2{\rm d}t\sim \\
& 2la\int_T^{[T]^1]}\left|\zf\right|^2{\rm d}t. 
\end{split}
\ee 
Consequently, in the case 
\be \label{5.2} 
l=l_1:\ 2l_1=\frac{1}{a}, 
\ee  
we obtain the next Theorem. 

\begin{mydef17}
On the Riemann hypothesis, the following asymptotic factorization formula 
\be \label{5.3} 
\begin{split}
& \int_T^{[T]^1]}\left|\zf\right|^2{\rm d}t\sim \\ 
& \int_{{N_1(T)}}^{N_1([T]^1)}\{-P(t)\}{\rm d}t\times \int_{[T]^1}^{[T+\frac{1}{a}]^1}\left|\zf\right|^2{\rm d}t,\ T\to\infty
\end{split}
\ee 
holds true.  
\end{mydef17} 

\subsection{} 

Let us put 
\be \label{5.4} 
Q_2(T)=\{(u,v)\in\mbb{R}^2:\ u\in(N_1(T),N_1([T]^1)),\ v\in ([T]^1,[T+\frac 1a]^1)\}. 
\ee 
Now, the formula (\ref{5.3}) implies by (\ref{2.6}), $r=1$, and (\ref{5.4}) the following result. 

\begin{mydef42}
On the Riemann hypothesis it is true that 
\be \label{5.5} 
\iint_{Q_2(T)}\{-P(u)\}\left|\zeta\left(\frac 12+iv\right)\right|^2{\rm d}u{\rm d}v=(1-c)T,\ T\to\infty. 
\ee 
\end{mydef42} 

\begin{remark}
Our asymptotic formula (\ref{5.5}) represents the result of a new kind for the interaction 
\be \label{5.6} 
-P(u) \ \leftrightarrow \ \left|\zeta\left(\frac 12+iv\right)\right|^2
\ee 
i.e. for the interaction between the remainder of the prime-number law and the quadrate of modulus of the Riemann zeta-function on the critical line. 
\end{remark} 

\subsection{} 

Next, if we put 
\be \label{5.7} 
\delta=\frac{1}{T\ln^kT};\ \{\ln^kT\}_{k=1}^\infty
\ee 
for every fixed $k\in\mbb{N}$ in (\ref{1.3}), then we obtain the following formula 
\be \label{5.8} 
\int_2^{N_2^k(T\ln^kT)}\{-P(t)\}{\rm d}t=\{1+o(1)\}aT\ln^kT,\ T\to\infty, 
\ee   
where 
\be \label{5.9} 
N_2^k(T\ln^kT)=N\left(\frac{1}{T\ln^kT}\right). 
\ee 
Now, we apply two formulae (\ref{3.11}) and (\ref{2.5}), $r=1$ in (\ref{5.8}) and we obtain the next result 
\be \label{5.10} 
\begin{split}
& \int_2^{N_2^k(T\ln^kT)}\{-P(t)\}{\rm d}t=\{1+o(1)\}\frac{a}{2l(1-c)}\times \\ 
& \int_T^{[T]^1}\left|\zf\right|^2{\rm d}t\times \int_{[T]^k}^{[T+2l]^k}\prod_{r=0}^{k-1}\left|\zfvpr\right|^2{\rm d}t,\ T\to\infty. 
\end{split}
\ee 
If we make the next choice for the value $2l$: 
\be \label{5.11} 
2l_2=\frac{a}{1-c}, 
\ee  
then we obtain the following theorem. 

\begin{mydef18}
On the Riemann hypothesis the following asymptotic factorization formula holds true for every fixed $k\in\mbb{N}$: 
\be \label{5.12} 
\begin{split}
& \int_2^{N_2^k(T\ln^kT)}\{-P(t)\}{\rm d}t\sim \\ 
& \int_T^{[T]^1}\left|\zf\right|^2{\rm d}t\times \int_{[T]^k}^{[T+\frac{a}{1-c}]^k}\prod_{r=0}^{k-1}\left|\zfvpr\right|^2{\rm d}t, \\ 
& T\to\infty, \ [G]^k=\vp_1^{-k}(G). 
\end{split}
\ee 
\end{mydef18} 

\begin{remark}
Let us remind explicitly that we have used only one of the sets of scales of infinity\footnote{See (\ref{5.7}).} 
\bdis 
\ln T,\ \ln^2T, \dots, \ \ln^kT,\ \dots 
\edis 
in the sense of Hardy's Pure Mathematics, see \cite{1}, pp. 346, 350, 370, that is sufficient for our purpose. 
\end{remark} 

\section{The asymptotic multiplicative formula as multiplicative puzzle on a set of integrals} 

\subsection{} 

Now we put in (\ref{1.3}) 
\be \label{6.1} 
\delta=\frac{1}{\ln^kT}
\ee  
for every fixed $k\in\mbb{N}$ and this gives to us the following formula 
\be \label{6.2} 
\int_2^{N_3^k(\ln^kT)}\{-P(t)\}{\rm d}t=\{1+o(1)\}\ln^kT, 
\ee  
where 
\be \label{6.3} 
N_3^k(\ln^kT)=N\left(\frac{1}{\ln^kT}\right). 
\ee 
Next, we use in (\ref{6.2}) the partition 
\be \label{6.4} 
k=1+1+\dots+1 
\ee 
i.e. 
\be \label{6.5} 
\ln^kT=\ln T\times \ln T \times \dots \times \ln T, 
\ee  
and, after this, we use consecutively $k$-fold of the formula (\ref{3.11}), $k=1$ and obtain the following result. 
\be \label{6.6} 
\begin{split}
& \int_2^{N_3^k(\ln^kT)}\{-P(t)\}{\rm d}t=\{1+o(1)\}\frac{a}{\prod_{r=1}^k(2l_r)}\times \\ 
& \prod_{r=1}^{k}\int_{[T]^1}^{[T+2l_r]^1}\left|\zf\right|^2{\rm d}t,\ l_1,\dots,l_k>0,\ T\to\infty 
\end{split}
\ee 
for every fixed $k\in\mbb{N}$. 

Of course, there is a continuum set of solutions 
\be \label{6.7} 
(2\bar{l}_1,2\bar{l}_2,\dots,2\bar{l}_k)\in (\mbb{R}^+)^k
\ee 
of the equation 
\be \label{6.8} 
\prod_{r=1}^k(2l_r)=a. 
\ee  
Namely, we can choose arbitrarily 
\be \label{6.9} 
(2\bar{l}_2,\dots,2\bar{l}_k)
\ee 
and then we have 
\be \label{6.10} 
2\bar{l}_1=\frac{a}{\prod_{r=2}^k(2l_r)}. 
\ee 
As a consequence we have the following theorem, see (\ref{6.3}). 

\begin{mydef19}
On the Riemann hypothesis it is true that: 
\begin{itemize}
	\item[(a)] for the integral 
	\be \label{6.11} 
	\int_2^{N_3^k(\ln^kT)}\{-P(t)\}{\rm d}t,\ k\in\mbb{N}, 
	\ee 
	\item[(b)] 
	for the continuum subset of vectors 
	\be \label{6.12} 
	\begin{split} 
	& (2\bar{l}_1,2\bar{l}_2,\dots,2\bar{l}_k)\in (\mbb{R}^+)^k:\ \prod_{r=1}^k(2l_r)=a, \\ 
	& a=e^{4.5},\ k\geq 2, 
	\end{split} 
	\ee 
	\item[(c)] 
	and for the set of integrals 
	\be \label{6.13} 
	\left\{\int_{[T]^1}^{[T+2\bar{l}_r]^1}\left|\zf\right|^2{\rm d}t\right\}_{r=1}^k
	\ee 
\end{itemize} 
there is a continuum set of asymptotic multiplicative formulas 
\be \label{6.14} 
\begin{split}
& \int_2^{N_3^k(\ln^kT)}\{-P(t)\}{\rm d}t\sim  \prod_{r=1}^k\int_{[T]^1}^{[T+2\bar{l}_r]^1}\left|\zf\right|^2{\rm d}t,\ T\to\infty
\end{split}
\ee 
for every fixed $k\geq 2$. 
\end{mydef19} 

\begin{remark}
The formula (\ref{6.14}) expresses the asymptotic solution of the problem to approximate the integral (\ref{6.11}) by means of the integrals of the class (\ref{6.13}). 
\end{remark} 

\begin{remark}
Asymptotic nonlinear composition (\ref{6.14}) of the integral (\ref{6.11}) by means of the integrals from the set (\ref{6.13}) can be viewed as a multiplicative puzzle on the set of integrals (\ref{6.13}).  
\end{remark} 

\subsection{} 

Here we give some examples of the formula\footnote{See (\ref{6.14}), $k=4$.} 
\be \label{6.15} 
\int_2^{N_3^4(\ln^4T)}\{-P(t)\}{\rm d}t\sim \prod_{r=1}^4\int_{[T]^1}^{[T+2\bar{l}_r]^1}\left|\zf\right|^2{\rm d}t,\ T\to\infty. 
\ee 
First, we have to get some solution 
\be \label{6.16} 
(2\bar{l}_1^s,\dots,2\bar{l}_4^s),\ s=1,\dots,s_0 
\ee 
of the equation\footnote{See (\ref{6.12}), $k=4$.} 
\be \label{6.17} 
\prod_{r=1}^4(2\bar{l}_r^s)=e^{4.5}. 
\ee 
We use the partitions $p(4)=5$ of the number $4$ i.e. 
\be \label{6.18} 
1+1+1+1,\ 2+1+1,\ 2+2,\ 3+1,\ 4 
\ee  
for this purpose. Corresponding solutions of the eq. (\ref{6.17}) read: 
\be \label{6.19} 
\begin{split}
& (2\bar{l}_1^1,\dots,2\bar{l}_4^1) = (e^0,e^1,e^{1.5},e^2), \\ 
& (2\bar{l}_1^2,\dots,2\bar{l}_4^2) = (e^{0.25},e^{0.25},e^1,e^3), \\ 
& (2\bar{l}_1^3,\dots,2\bar{l}_4^3) = (e^{0.25},e^{0.25},e^2,e^2), \\ 
& (2\bar{l}_1^4,\dots,2\bar{l}_4^4) = (e^{1.5},e^{1.5},e^{1.5},e^0), \\ 
& (2\bar{l}_1^5,\dots,2\bar{l}_4^5) = (e^{1.125},e^{1.125},e^{1.125},e^{1.125}), \\
\end{split}
\ee 
and, finally, corresponding nonlinear compositions (multiplicative puzzles) for the integral (\ref{6.15}) read:  
\bdis 
\int_2^{N_3^4(\ln^4T)}\{-P(t)\}{\rm d}t\sim 
\edis 
\be \label{6.20} 
\begin{split}
& \sim \int_{[T]^1}^{[T+e^0]^1}\left|\zf\right|^2{\rm d}t\times \int_{[T]^1}^{[T+e]^1}\left|\zf\right|^2{\rm d}t\times \\ 
& \int_{[T]^1}^{[T+e^{1.5}]^1}\left|\zf\right|^2{\rm d}t\times \int_{[T]^1}^{[T+e^2]^1}\left|\zf\right|^2{\rm d}t, 
\end{split}
\ee  
\be \label{6.21} 
\begin{split}
	& \sim \left\{\int_{[T]^1}^{[T+e^{0.25}]^1}\left|\zf\right|^2{\rm d}t\right\}^2\times \int_{[T]^1}^{[T+e]^1}\left|\zf\right|^2{\rm d}t\times \\ 
	& \int_{[T]^1}^{[T+e^{3}]^1}\left|\zf\right|^2{\rm d}t, 
\end{split}
\ee 
\be \label{6.22} 
\begin{split}
	& \sim \left\{\int_{[T]^1}^{[T+e^{0.25}]^1}\left|\zf\right|^2{\rm d}t\right\}^2\times \left\{\int_{[T]^1}^{[T+e^2]^1}\left|\zf\right|^2{\rm d}t\right\}^2,  
\end{split}
\ee 
\be \label{6.23} 
\begin{split}
	& \sim \left\{\int_{[T]^1}^{[T+e^{1.5}]^1}\left|\zf\right|^2{\rm d}t\right\}^3\times \int_{[T]^1}^{[T+e^0]^1}\left|\zf\right|^2{\rm d}t,  
\end{split}
\ee 
\be \label{6.24} 
\begin{split}
	& \sim \left\{\int_{[T]^1}^{[T+e^{1.125}]^1}\left|\zf\right|^2{\rm d}t\right\}^4,  \ T\to\infty. 
\end{split}
\ee 

\begin{remark} 
In our simple example we have used all the partitions of the number $4$. For bigger $k$ this method is no more applicable since $p(k)$ grows rapidly. For example 
\bdis 
p(200)= 3\ 372 \ 999 \ 029 \ 388. 
\edis 
\end{remark}

I would like to thank Michal Demetrian for his moral support of my study of Jacob's ladders.

\end{document}